\documentclass[oneside,british]{amsart}
\usepackage[T1]{fontenc}
\usepackage[latin1]{inputenc}
\usepackage{amssymb}

\makeatletter

\providecommand{\tabularnewline}{\\}

 \theoremstyle{plain}    
 \newtheorem{thm}{Theorem}[section]
 \numberwithin{equation}{section} 
 \numberwithin{figure}{section} 
 \theoremstyle{plain}    
 \newtheorem*{conjecture*}{Conjecture} 
 \theoremstyle{plain}    
 \newtheorem{prop}[thm]{Proposition} 
 \theoremstyle{definition}
 \newtheorem*{defn*}{Definition}
 \usepackage{verbatim}
 \theoremstyle{plain}    
 \newtheorem{lem}[thm]{Lemma} 
 \theoremstyle{definition}
 \newtheorem{defn}[thm]{Definition}

\usepackage{alex}
\renewcommand{\theenumi}{\alph{enumi}}

\usepackage{babel}
\makeatother
\begin{document}

\title{Representations of Reductive Groups over Quotients of Local Rings }

\author{Alexander Stasinski}

\address{School of Mathematical Sciences\\
University of Nottingham\\
NG7 2RD\\
England}

\email{pmxas@maths.nott.ac.uk}

\begin{abstract}
In some recent work, Lusztig outlined a generalisation of the construction
of Deligne and Lusztig to reductive groups over finite rings coming
from the ring of integers in a local field, modulo some power of the
maximal ideal. Lusztig conjectures that all irreducible representations
of these groups are contained in the cohomology of a certain family
of varieties. We show that, contrary to what was expected, there exist
representations that cannot be realised by the varieties given by
Lusztig. Moreover, we show how the remaining representations in the
case under consideration can be realised in the cohomology of a different
kind of variety. This may suggest a way to reformulate Lusztig's conjecture.
\end{abstract}
\maketitle

\section{Introduction}

Let $G$ be a connected reductive group, defined over a finite field
$\mathbb{F}_{q}$ of characteristic $p$. The celebrated work of Deligne
and Lusztig (\cite{delignelusztig}) uses methods of $l$-adic \'etale
cohomology to show that many complex irreducible representations of
the finite group $G(\mathbb{F}_{q})$ are parametrised by characters
of maximal tori in $G(\mathbb{F}_{q})$, and moreover, that every
irreducible representation appears in the cohomology of a certain
kind of variety.

Now let $K$ be a local field (of any characteristic) with finite
residue field $\mathbb{F}_{q}$, and ring of integers $\mathcal{O}_{K}$.
Fix an algebraic closure of $K$, and let $\mathcal{O}$ be the ring
of integers of the maximal unramified extension $K^{\text{ur}}$ of
$K$ with residue field $\mathbb{F}$, an algebraic closure of $\mathbb{F}_{q}$.
Denote by $\epsilon$ a fixed prime element in $\mathcal{O}_{K}$.
It is also a prime element in $\mathcal{O}$.

Assume that $X$ is an affine variety over $\mathbb{F}$, and let
$r\geq1$ be an integer. We set\[
X_{r}=X(\mathcal{O}/(\epsilon^{r})).\]
Thus, if $X$ is the common zeroes of the polynomials $f_{i}(x_{1},\dots,x_{n})$,
$i=1,\dots,m$, then $X_{r}$ is the set of all $(a_{1},\dots,a_{n})\in(\mathcal{O}/(\epsilon^{r}))^{n}$
such that $f_{i}(a_{1},\dots,a_{n})=0$ for $i=1,\dots,m$. This makes
sense, since $\mathcal{O}/(\epsilon^{r})$ is an $\mathbb{F}$-algebra.
Then $X\mapsto X_{r}$ is a functor from the category of affine varieties
over $\mathbb{F}$ into itself. 

For any $r\geq r'\geq0$ the reduction map $\mathcal{O}/(\epsilon^{r})\rightarrow\mathcal{O}/(\epsilon^{r'})$
induces a morphism $\rho_{r,r'}:X_{r}\rightarrow X_{r'}$. We will
denote the map $\rho_{r,1}$ by $\rho_{r}$. Consider the reductive
algebraic group $G$ over $\mathbb{F}$. Then $G_{r}$ is naturally
an algebraic group over $\mathbb{F}$, and $\rho_{r,r'}:G_{r}\rightarrow G_{r'}$
is a surjective homomorphism of algebraic groups with kernel $G_{r}^{r'}$.
Thus we have an exact sequence\[
1\longrightarrow G_{r}^{r'}\longrightarrow G_{r}\xrightarrow{\rho_{r,r'}}G_{r'}\longrightarrow1.\]
The injection $\mathcal{O}/(\epsilon^{r'})\rightarrow\mathcal{O}/(\epsilon^{r})$
induces a function $i_{r',r}:G_{r'}\rightarrow G_{r}$ such that $\rho_{r,r'}\circ i_{r',r}$
is the identity map on $G_{r'}$. In the case where $r'=1$ and $K$
is a local field of positive characteristic there is an inclusion
of $\mathbb{F}$-algebras $\mathcal{O}/(\epsilon)\rightarrow\mathcal{O}/(\epsilon^{r})$,
and $i_{1,r}$ is an injective homomorphism, so that the above exact
sequence splits. When $K$ is of characteristic zero, the exact sequence
is just a group extension. We will denote $i_{1,r}$ by $i_{r}$.

Let $F:G\rightarrow G$ be the Frobenius morphism corresponding to
the $\mathbb{F}_{q}$-rational structure of $G$. The map $F$ induces
a homomorphism $F:G_{r}\rightarrow G_{r}$ which is the Frobenius
map for an $\mathbb{F}_{q}$-rational structure on $G_{r}$. Let $T$
be a maximal torus in $G$ contained in an $F$-stable Borel subgroup
$B$, with unipotent radical $U$. Fix $r\geq1$ as above, and let
$x\in G_{r}$. Consider the affine variety\[
X_{x}=\{ g\in G_{r}\mid g^{-1}F(g)\in xU_{r}\}.\]
The group $G_{r}^{F}$ of fixed points under the Frobenius map can
be identified with the $\mathbb{F}_{q}$-points of $G_{r}$, and is
thus a finite group. The group $G_{r}^{F}$ acts on $X_{x}$ by left
multiplication, and thus by functoriality, the $l$-adic cohomology
with compact support $H_{c}^{i}(X_{x},\overline{\mathbb{Q}}_{l})$
for $l\neq p$, has the structure of a complex representation of $G_{r}^{F}$.
Note that $\overline{\mathbb{Q}}_{l}\simeq\mathbb{C}$ (non-canonically). 

The above construction is an extension of the construction of Deligne
and Lusztig (which is the case $r=1$), and was first mentioned by
Lusztig in \cite{lusztig}, and then recently developed further in
\cite{lusztigpreprint}. This construction is the first natural step
to an extension of the Deligne-Lusztig construction to reductive groups
over local fields, a problem which has important arithmetic implications.

In the paper \cite{lusztigpreprint}, Lusztig proves an orthogonality
formula for certain virtual representations of the group $G_{r}^{F}$,
in the case where $K$ is a local field of positive characteristic.
This result is an extension of an important result of Deligne and
Lusztig in the case $r=1$, and was anticipated in \cite{lusztig}. 

It was shown in \cite{delignelusztig} that for the case $r=1$ every
irreducible representation of $G^{F}$ appears in the cohomology of
some variety $X_{w}$, where $w$ is an element of the Weyl group
of $G$ (the varieties are independent of the lift of $w$ to an element
of $G$). As pointed out in \cite{lusztigpreprint}, this is no longer
the case for $r\geq2$. In the end of \cite{lusztigpreprint}, Lusztig
states the following

\begin{conjecture*}
Any irreducible representation of $G_{r}^{F}$ appears in the virtual
representation $\sum_{i\geq0}(-1)^{i}H_{c}^{i}(X_{x},\overline{\mathbb{Q}}_{l})$
for some $x\in G_{r}$.
\end{conjecture*}
In the following we will show that the conjecture does not hold for
$G=\textrm{SL}_{2}$, with $K$ a local field of positive characteristic
with $q$ odd, and $r=2$, i.e., for the group $G_{2}^{F}=\textrm{SL}_{2}(\mathbb{F}_{q}[[\epsilon]]/(\epsilon^{2}))$.
In the end of this paper we will show that the missing irreducible
representations of this group are realised in the cohomology of a
certain variety, not of the form $X_{x}$ for any $x\in G_{2}$. These
representations are parametrised by certain characters of two different
subgroups of $G_{2}^{F}$.

So far nothing seems to be known about Lusztig's conjecture in the
case where $K$ is a local field of mixed characteristic. This is
the focus of work in progress (\cite{alex}).

\section{\label{sec:The-Reps-of}The Representations of $\textrm{SL}_{2}(\mathbb{F}_{q}[[\epsilon]]/(\epsilon^{2}))$}

Recall the notation of the previous section. We specialise the discussion
to the case $G=\textrm{SL}_{2}$, $K$ of positive characteristic
with $q$ odd, and $r=2$. The irreducible representations of $G_{2}^{F}$
can be classified using the fact that $G_{2}^{F}$ is a semidirect
product of $G^{F}$, and a group $N$ isomorphic to $(\mathbb{F}_{q}^{+})^{3}$
(three copies of the additive group of the finite field). The following
table of representations of $G_{2}^{F}$ was given in \cite{lusztigpreprint}.
The first column indicates the dimension, and the second column indicates
the number of representations of that dimension.

\begin{center}\bigskip\begin{tabular}{c|c}
dim&
\#\tabularnewline
\hline 
$1$&
$1$\tabularnewline
$q$&
$1$\tabularnewline
$q+1$&
$(q-3)/2$\tabularnewline
$(q+1)/2$&
$2$\tabularnewline
$q-1$&
$(q-1)/2$\tabularnewline
$(q-1)/2$&
$2$\tabularnewline
$q^{2}+q$&
$(q-1)^{2}/2$\tabularnewline
$q^{2}-q$&
$(q^{2}-1)/2$\tabularnewline
$(q^{2}-1)/2$&
$2q$\tabularnewline
\end{tabular}\bigskip\end{center}

A calculation using the fact that the sum of the squares of the dimensions
of the irreducible representations of a finite group equals the order
of the group, shows that this table is incorrect. We will show in
the next section that the correct value for the last entry in the
second column is $4q$. We therefore have to consider $2q$ representations
in addition to those considered by Lusztig.

Consider $G_{2}^{F}$ as a semidirect product of $G^{F}$ and $N$,
where $N$ is the normal abelian subgroup of $G_{2}^{F}$ consisting
of matrices of the form\[
\Bigl{\{}g=\begin{pmatrix}1+a\epsilon & b\epsilon\\
c\epsilon & 1+d\epsilon\end{pmatrix}\mid\det g=1\Bigr{\}}=\Bigl{\{}\begin{pmatrix}1+a\epsilon & b\epsilon\\
c\epsilon & 1-a\epsilon\end{pmatrix}\Bigr{\}}.\]
 The method of describing the representations of a semidirect product
by an abelian group is given in \cite{serrelinrep}. We give here
a summary of this method for the case under consideration.

Let $X=\Hom(N,\mathbb{C}^{\times})$. The group $G_{2}^{F}$ acts
on $X$ by\[
(s\chi)(n)=\chi(s^{-1}ns)\quad\textrm{for }s\in G_{2}^{F},\chi\in X,n\in N.\]
Let $(\chi_{i})_{i\in X/G^{F}}$ be a system of representatives for
the orbits of $G^{F}$ in $X$. For each $i\in X/G^{F}$, let $(G^{F})_{i}=\textrm{Stab}_{G^{F}}(\chi_{i})$
and let $(G_{2}^{F})_{i}=(G^{F})_{i}\cdot N$ be the corresponding
subgroup of $G_{2}^{F}$. The character $\chi_{i}$ can be extended
to $(G_{2}^{F})_{i}$ by setting\[
\chi_{i}(gn)=\chi_{i}(n)\quad\textrm{for }g\in(G^{F})_{i},n\in N.\]
Now let $\rho$ be an irreducible representation of $(G^{F})_{i}$.
By composing $\rho$ with the canonical projection $(G_{2}^{F})_{i}\rightarrow(G^{F})_{i}$
we obtain an irreducible representation $\tilde{\rho}$ of $(G_{2}^{F})_{i}$.
Finally, by taking the tensor product of $\tilde{\rho}$ and $\chi_{i}$
we obtain an irreducible representation $\tilde{\rho}\otimes\chi_{i}$
of $(G_{2}^{F})_{i}$. Let $\theta_{i,\rho}$ be the corresponding
induced representation of $G_{2}^{F}$. Then we have the following
result (cf. \cite{serrelinrep}, Proposition 25, 8.2.)

\begin{prop}
~
\begin{enumerate}
\item $\theta_{i,\rho}$ is irreducible.
\item If $\theta_{i,\rho}$ and $\theta_{i',\rho'}$ are isomorphic, then
$i=i'$ and $\rho\simeq\rho'$.
\item Every irreducible representation of $G_{2}^{F}$ is isomorphic to
one of the $\theta_{i,\rho}$.
\end{enumerate}
\end{prop}
We now show how to classify all the irreducible representations of
$G_{2}^{F}$ of degree $(q^{2}-1)/2$ using the above method. Let
$\psi:\mathbb{F}_{q}^{+}\rightarrow\mathbb{C}^{\times}$ be a nontrivial
additive character, and consider the character\[
\chi_{\psi}:N\longrightarrow\mathbb{C}^{\times},\quad\chi_{\psi}\begin{pmatrix}1+a\epsilon & b\epsilon\\
c\epsilon & 1-a\epsilon\end{pmatrix}=\psi(c).\]
Computing the stabiliser of $\chi_{\psi}$ in $G^{F}$, we get\begin{eqnarray*}
\textrm{Stab}_{G^{F}}(\chi_{\psi}) & = & \Bigl{\{}\begin{pmatrix}x & y\\
z & w\end{pmatrix}\in G^{F}\mid\chi_{\psi}\begin{bmatrix}\begin{pmatrix}x & y\\
z & w\end{pmatrix}^{-1}\begin{pmatrix}1+a\epsilon & b\epsilon\\
c\epsilon & 1-a\epsilon\end{pmatrix}\begin{pmatrix}x & y\\
z & w\end{pmatrix}\end{bmatrix}\\
 & = & \chi_{\psi}\begin{pmatrix}1+a\epsilon & b\epsilon\\
c\epsilon & 1-a\epsilon\end{pmatrix},\forall a,b,c\in\mathbb{F}_{q}^{+}\Bigr{\}}\\
 & = & \Bigl{\{}\begin{pmatrix}x & y\\
z & w\end{pmatrix}\in G^{F}\mid\psi(x^{2}c-z^{2}b+2xza)=\psi(c),\forall a,b,c\in\mathbb{F}_{q}^{+}\Bigr{\}}\\
 & = & \Bigl{\{}\begin{pmatrix}x & y\\
0 & x^{-1}\end{pmatrix}\in G^{F}\mid x^{2}=1\Bigr{\}}=\Bigl{\{}\begin{pmatrix}\pm1 & y\\
0 & \pm1\end{pmatrix}\mid y\in\mathbb{F}_{q}\Bigr{\}}.\end{eqnarray*}
Since $\textrm{Stab}_{G^{F}}(\chi_{\psi})$ is an abelian group of
order $2q$, we obtain $2q$ irreducible representations of $G_{2}^{F}$
of dimension $(q^{2}-1)/2$ by inducing certain $1$-dimensional representations
of $\textrm{Stab}_{G^{F}}(\chi_{\psi})\cdot N$, according to the
method described above. Now the problem is to determine the orbits
of all characters of the form $\chi_{\psi}$ under the action of $G^{F}$.

Let $\zeta=e^{2\pi i/p}$ and let $\Tr:\mathbb{F}_{q}\rightarrow\mathbb{F}_{p}$
be the absolute trace. It is well known that every character $\mathbb{F}_{q}^{+}\rightarrow\mathbb{C}^{\times}$
is of the form\[
\psi_{a}(x)=\zeta^{\Tr(ax)},\]
for some $a\in\mathbb{F}_{q}$, and that $\psi_{a}\neq\psi_{b}$ for
$a\neq b$.

We now consider the action of $G^{F}$ on the characters $\psi_{a}$,
where $a\in\mathbb{F}_{q}^{\times}$. For $k,m\in\mathbb{F}_{q}^{\times}$
we see that if $\chi_{\psi_{k}}$ lies in the same orbit as $\chi_{\psi_{m}}$,
then there exists an element $g=\left(\begin{smallmatrix}x & y\\
z & w\end{smallmatrix}\right)\in G^{F}$such that\[
(g\chi_{\psi_{k}})\begin{pmatrix}1+a\epsilon & b\epsilon\\
c\epsilon & 1-a\epsilon\end{pmatrix}=\chi_{\psi_{m}}\begin{pmatrix}1+a\epsilon & b\epsilon\\
c\epsilon & 1-a\epsilon\end{pmatrix}\quad\textrm{for all }a,b,c\in\mathbb{F}_{q}.\]
This in turn implies that\[
\psi_{k}(x^{2}c-z^{2}b+2xza)=\psi_{m}(c)\quad\textrm{for all }a,b,c\in\mathbb{F}_{q},\]
and so\[
\psi_{k}(x^{2}c)=\psi_{m}(c)\quad\textrm{for all }c\in\mathbb{F}_{q}.\]
Thus we must have $kx^{2}=m$. Conversely, if $kx^{2}=m$ holds, then
clearly $\chi_{\psi_{k}}$ and $\chi_{\psi_{m}}$ lie in the same
$G^{F}$-orbit. 

We can now identify two orbits. We will denote the subgroup of squares
in $\mathbb{F}_{q}^{\times}$ by $\mathbb{F}_{q}^{\times2}$. First,
assume that $k\in\mathbb{F}_{q}^{\times2}$. Then $kx^{2}$ is a square,
and any square in $\mathbb{F}_{q}^{\times}$ can be expressed in this
way for some $x\in\mathbb{F}_{q}^{\times}$, since $\mathbb{F}_{q}^{\times2}$
is a group. Thus $\{\chi_{\psi_{k}}\mid k\in\mathbb{F}_{q}^{\times2}\}$
is an orbit. On the other hand, if $k$ is not a square, then neither
is $kx^{2}$, and all nonsquares in $\mathbb{F}_{q}^{\times}$ can
be expressed in this way for some $x\in\mathbb{F}_{q}^{\times}$,
since $kx^{2}=kx'^{2}\Rightarrow x^{2}=x'^{2}$, and so we get all
$|\mathbb{F}_{q}^{\times2}|$ distinct nonsquares. These elementary
considerations can be summed up as \[
G^{F}\backslash\{\chi_{\psi_{k}}\mid k\in\mathbb{F}_{q}^{\times}\}\simeq\mathbb{F}_{q}^{\times}/\mathbb{F}_{q}^{\times2}\simeq\{\pm1\}.\]

The above discussion shows that there are $4q$ irreducible representations
of dimension $(q^{2}-1)/2$. Half of them correspond to the orbit
$\{\chi_{\psi_{k}}\mid k\in\mathbb{F}_{q}^{\times2}\}$, which is
the one containing $\chi_{\psi_{1}}$, and the other half correspond
to the orbit corresponding to nonsquares in $\mathbb{F}_{q}^{\times}$.

\section{Some Lemmas}

In this section we fix some notation and collect some results from
Deligne-Lusztig theory. Most results from the finite field situation
($r=1$) hold for all $r$. All results in this section except for
the last one, were proved for $r=1$ by Deligne and Lusztig (cf. \cite{delignelusztig},
and \cite{lusztigexpo}). 

Let $L:G_{r}\rightarrow G_{r}$, $L(g)=g^{-1}F(g)$ denote the Lang
map. We have \[
L^{-1}(xU_{r})=\{ g\in G_{r}\mid g^{-1}F(g)\in xU_{r}\}=X_{x}.\]
From now on we will denote $H_{c}^{i}(X,\overline{\mathbb{Q}}_{l})$
simply by $H_{c}^{i}(X)$.

\begin{defn*}
Let $G$ be a finite group that acts on the varieties $X$ and $Y$.
We write $X\sim Y$ if any irreducible representation of $G$ appears
in $\sum_{i\geq0}(-1)^{i}H_{c}^{i}(X)$ if and only if it appears
in $\sum_{i\geq0}(-1)^{i}H_{c}^{i}(Y)$.
\end{defn*}
Note that the relation just defined is an equivalence relation.

\begin{lem}
\label{lem:X/U}Suppose that $f:X\rightarrow Y$ is a morphism of
varieties such that for some $m\geq0$, the fibre $f^{-1}(y)$ is
isomorphic to affine $m$-space, for all $y\in Y$. Let $g,g'$ be
automorphisms of finite order of $X,Y$ such that $fg=g'f$. Then
$X\sim Y$.
\end{lem}
\begin{proof}
See \cite{lusztigexpo}, Lemma 1.9.
\end{proof}
As pointed out in \cite{lusztigpreprint}, for arbitrary $r$ it is
enough for representation theoretic purposes to consider varieties
$X_{x}$ where $x$ runs through a set of double coset representatives
$U_{r}\backslash G_{r}/U_{r}$. This follows from Lemma \ref{lem:X/U}
and the following result.

\begin{lem}
\label{lem:UxU}The inclusion $L^{-1}(xU_{r})\hookrightarrow L^{-1}(U_{r}xU_{r})$
induces an isomorphism\[
L^{-1}(xU_{r})/U_{r}\cap xU_{r}x^{-1}\longiso L^{-1}(U_{r}xU_{r})/U_{r},\]
commuting with the action of $G_{r}^{F}$ on both varieties.
\end{lem}
\begin{proof}
Denote by $f$ the composition of the maps $L^{-1}(xU_{r})\hookrightarrow L^{-1}(U_{r}xU_{r})\rightarrow L^{-1}(U_{r}xU_{r})/U_{r}$,
where the latter is the natural projection. Clearly $f$ is surjective,
because if $gU_{r}\in L^{-1}(U_{r}xU_{r})/U_{r}$, with $L(g)\in uxu'$
for $u,u'\in U_{r}$, then $L(gu)=u^{-1}uxu'F(u)\in xU_{r}$, so $gu\in L^{-1}(xU_{r})$,
and $f(gu)=gU_{r}$.

On the other hand, the fibre of $f$ at $gU_{r}$ is equal to $\{ gv\in L^{-1}(xU_{r})\mid v\in U_{r}\}=\{ gv\mid v^{-1}L(g)F(v)\in xU_{r},\  v\in U_{r}\}=\{ gv\mid v^{-1}ux\in xU_{r},v\in U_{r}\}=\{ gv\mid v^{-1}u\in U_{r}\cap xU_{r}x^{-1}\}=\{ gv\mid v=u\bmod U_{r}\cap xU_{r}x^{-1}\}$.
Factoring $L^{-1}(xU_{r})$ by $U_{r}\cap xU_{r}x^{-1}$ therefore
gives an isomorphism, which commutes with the action of $G_{r}^{F}$,
by the naturality of $f$.
\end{proof}
\begin{lem}
\label{lem:F(x)UF(x)^-1}Let $x\in G_{r}$ be an arbitrary element,
and let $\lambda$ be an element such that $L(\lambda)=x$. Then there
is an isomorphism\[
L^{-1}(xU_{r})\longiso L^{-1}(F(\lambda)U_{r}F(\lambda)^{-1}),\quad g\longmapsto g\lambda^{-1},\]
commuting with the action of $G_{r}^{F}$. 
\end{lem}
\begin{proof}
Let $g\in L^{-1}(xU_{r})$. Then $L(g\lambda^{-1})=\lambda L(g)F(\lambda)^{-1}\in\lambda xU_{r}F(\lambda)^{-1}=F(\lambda)x^{-1}xU_{r}F(\lambda)^{-1}=F(\lambda)U_{r}F(\lambda)^{-1}$.
It is clear that this map is a morphism of varieties, and it has an
obvious inverse.
\end{proof}
In the case $r=1$, the Bruhat decomposition says that $B\backslash G/B$
is in bijection with the Weyl group, and that a set of double coset
representatives can be taken in $N_{G}(T)$, the normaliser of $T$.
Suppose that $utwt'u'$ is an arbitrary element in $BwB$, for $w\in N_{G}(T)$.
Then\[
L^{-1}(Uutwt'u'U)=L^{-1}(Utwt'U)\sim L^{-1}(twt'U)=L^{-1}(t''wU),\]
for some $t''\in T$. Since $t\mapsto wF(t)w^{-1}$ is a Frobenius
map on $T$, Lang's theorem says that there exists $\lambda\in T$
such that $\lambda^{-1}wF(\lambda)w^{-1}=t''$. The map\[
L^{-1}(t''wU)\longrightarrow L^{-1}(wU),\quad g\longmapsto g\lambda^{-1}\]
is then an isomorphism of varieties, preserving the action of $G^{F}$.
Indeed, if $g$ is an element in $L^{-1}(t''wU)$, then $L(g\lambda^{-1})\in\lambda t''wU=wF(\lambda)w^{-1}t''^{-1}t''wUF(\lambda)^{-1}=wF(\lambda)UF(\lambda)=wU$. 

Because of this, in the case $r=1$ it is enough to consider varieties
attached to elements of the Weyl group, among all varieties $X_{x}$,
$x\in G$. This is no longer enough for $r>1$.

\begin{defn}
Let $G$ be a finite group that acts on the left on the variety $X$,
and let $A$ be a finite abelian group that acts on the right on $X$.
For any character $\theta:A\rightarrow\overline{\mathbb{Q}}_{l}^{\times}$
and $g\in G$, we define a virtual character of $G$ by\[
\mathscr{L}(g,X)_{\theta}=\sum_{i\geq0}(-1)^{i}\Tr(g,H_{c}^{i}(X)_{\theta}).\]

\end{defn}
Here we use the notation $V_{\theta}=\{ v\in V\mid va=\theta(a)v,\textrm{ for all }a\in A\}$,
for $V$ a finite dimensional right $\overline{\mathbb{Q}}_{l}[A]$-module.
Since for any such $V$ we have $V=\bigoplus_{\theta}V_{\theta}$,
we have a virtual character\[
\mathscr{L}(g,X)=\bigoplus_{\theta}\mathscr{L}(g,X)_{\theta}=\sum_{i\geq0}(-1)^{i}\Tr(g,H_{c}^{i}(X)),\]
classically called the \emph{Lefschetz number}. We will make use of
the following results.

\begin{lem}
\label{lem:char-formula}Let $G$, $A$, and $X$ be as in the above
definition. Then for any $g\in G$ we have \[
\mathscr{L}(g,X)_{\theta}=\frac{1}{|A|}\sum_{a\in A}\theta(a^{-1})\mathscr{L}((g,a),X),\]
where $(g,a)$ acts on $X$ by $x\mapsto gxa$.
\end{lem}
\begin{proof}
See \cite{carter}, Proposition 7.2.3.
\end{proof}
\begin{lem}
\label{lem:fixed-point}Let $G$ and $X$ be as above, and assume
that $X$ is a finite set. Then $H_{c}^{i}(X)=0$ if $i\neq0$, and
$H_{c}^{0}(X)\simeq\overline{\mathbb{Q}}_{l}[X]$ is a permutation
representation of $G$. Furthermore the character of this representation
is given by $\mathscr{L}(g,X)=|X^{g}|$.
\end{lem}
\begin{proof}
See e.g. \cite{dignemichel}, Proposition 10.8.
\end{proof}
\begin{lem}
\label{lem:partition}Let $G$ be a finite group acting on the variety
$X$, and let $X=\coprod_{i\in I}X_{i}$ be a finite partition of
$X$ into disjoint, closed subsets. Assume that $G$ permutes the
subsets $X_{i}$ among them in such a way that the action of $G$
on $I$ is transitive. Let $H=\Stab_{G}(X_{i_{0}})=\{ g\in G\mid gX_{i_{0}}=X_{i_{0}}\}$
for a fixed $i_{0}\in I$. Then the generalised character $g\mapsto\mathscr{L}(g,X)$
of $G$ is induced by the generalised character $h\mapsto\mathscr{L}(h,X_{i_{0}})$
of $H$.
\end{lem}
\begin{proof}
See \cite{lusztigexpo}, Lemma 1.7.
\end{proof}
We end with a new result which has no nontrivial analogue for $r=1$.

\begin{lem}
Let $r>1$, $x\in G_{r}$, and consider the variety $L^{-1}(xU_{r})$.
Then the projection map $\rho_{r}:G_{r}\rightarrow G$ induces an
isomorphism \[
(G_{r}^{1})^{F}\backslash L^{-1}(xU_{r})\longiso L^{-1}(\rho_{r}(x)U),\]
commuting with the action of $G_{2}^{F}$. 
\end{lem}
\begin{proof}
Let $f$ be the map $L^{-1}(xU_{r})\rightarrow L^{-1}(\rho_{r}(x)U)$,
given by $g\mapsto\rho_{r}(g)$. This map is surjective because if
$g\in L^{-1}(\rho_{r}(x)U)$, then $i_{r}(g)\in i_{r}(L^{-1}(\rho_{r}(x)U))\subset L^{-1}(xU_{r})$,
and $f(i_{r}(g))=g$. The fibre of $f$ at $g$ is equal to \[
\{ a\cdot i_{r}(g)\in L^{-1}(xU_{r})\mid\rho_{r}(a)=1\}=\{ a\cdot i_{r}(g)\mid a\in(G_{r}^{1})^{F}\},\]
and this shows that factoring $L^{-1}(xU_{r})$ by $(G_{r}^{1})^{F}$
gives an isomorphism which commutes with the action of $G_{r}^{F}$,
by the naturality of $f$.
\end{proof}

\section{Lusztig's Conjecture for $\textrm{SL}_{2}(\mathbb{F}_{q}[[\epsilon]]/(\epsilon^{2}))$}

In the rest of this paper we will focus on the case $G=\textrm{SL}_{2}$,
$K$ of positive characteristic with $q$ odd, and $r=2$. In the
following we will show that only $4$ of the representations of dimension
$(q^{2}-1)/2$ of $G_{2}^{F}$ can be realised by varieties of the
form $X_{x}$ for $x\in G_{2}$.

By Lemma \ref{lem:UxU}, all representations that can be realised
by varieties $X_{x}$, $x\in G_{2}$, can also be realised by varieties
$X_{x}$, where $x$ is a representative in $U_{2}\backslash G_{2}/U_{2}$.
Therefore, the first place to start looking for representations is
in the cohomology of varieties corresponding to representatives of
$B_{2}\backslash G_{2}/B_{2}$. Such a set of representatives is given
by\[
\Bigl{\{}w_{1}=\begin{pmatrix}1 & 0\\
0 & 1\end{pmatrix},\  w_{2}=\begin{pmatrix}0 & -1\\
1 & 0\end{pmatrix},\  e=\begin{pmatrix}1 & 0\\
\epsilon & 1\end{pmatrix}\Bigr\}.\]

\begin{prop}
Any irreducible representation that can be realised by a variety $X_{x}$
for $x\in G_{2}$, can also be realised by a variety $X_{x}$, where
either $x=w_{1},w_{2}$, or $x=\left(\begin{smallmatrix}1 & 0\\
k\epsilon & 1\end{smallmatrix}\right)$, for some $k\in\mathbb{F}^{\times}$.
\end{prop}
\begin{proof}
Every element $x\in G_{2}$ lies in a $B_{2}$--$B_{2}$ double coset
corresponding to one of the elements $w_{1},w_{1},e$ above. The elements
$w_{1}$ and $w_{2}$ normalise $T_{2}$, so as in the case of $r=1$,
for any element $x\in B_{2}w_{1}B_{2}=B_{2}$ we have $X_{x}\sim X_{w_{1}}$,
and for any $y\in B_{2}w_{2}B_{2}$ we have $X_{y}\sim X_{w_{2}}$. 

In contrast, the element $e$ does not normalise $T_{2}$, so we cannot
a priori draw the same conclusions as above. Assume that $x=utet'u'$,
where $u,u'\in U_{2}$ and $t,t'\in T_{2}$. Then $L^{-1}(utet'u'U_{2})\sim L^{-1}(U_{2}tet'U_{2})\sim L^{-1}(tet'U_{2})$,
and by Lemma \ref{lem:F(x)UF(x)^-1} we have $L^{-1}(tet'U_{2})\sim L^{-1}(F(\lambda)U_{2}F(\lambda)^{-1})$,
where $L(\lambda)=tet'$. We can assume that $\lambda$ has the form
\[
\lambda=\begin{pmatrix}t_{0}+t_{1}\epsilon & 0\\
u\epsilon & (t_{0}+t_{1}\epsilon)^{-1}\end{pmatrix}\quad\textrm{for some }u,t_{0}\in\mathbb{F}_{q}^{\times},t_{1}\in\mathbb{F}_{q}.\]
Since we can write $\lambda=e't''$, where $e'=\left(\begin{smallmatrix}1 & 0\\
ut_{0}^{-1}\epsilon & 1\end{smallmatrix}\right)$, and $t''=\left(\begin{smallmatrix}t_{0}+t_{1}\epsilon & 0\\
0 & (t_{0}+t_{1}\epsilon)^{-1}\end{smallmatrix}\right)$, we get $L^{-1}(F(\lambda)U_{2}F(\lambda)^{-1})=L^{-1}(F(e't'')U_{2}F(e't'')^{-1})=L^{-1}(F(e')U_{2}F(e')^{-1})\sim L^{-1}(L(e')U_{2})$.
The element $L(e')=(e')^{-1}F(e')$ is obviously of the form $\left(\begin{smallmatrix}1 & 0\\
k\epsilon & 1\end{smallmatrix}\right)$, for some $k\in\mathbb{F}^{\times}$. Thus, for every $x\in B_{2}eB_{2}$
we have $X_{x}\sim L^{-1}(\left(\begin{smallmatrix}1 & 0\\
k\epsilon & 1\end{smallmatrix}\right)U_{2})$, for some $k\in\mathbb{F}^{\times}$, and any such $k$ appears for
some $x$.
\end{proof}
It is clear that the varieties corresponding to elements $\left(\begin{smallmatrix}1 & 0\\
k\epsilon & 1\end{smallmatrix}\right)$, $k\in\mathbb{F}^{\times}$, are not all essentially different. Namely,
if $\kappa\in G_{2}$ is such that $L(\kappa)=\left(\begin{smallmatrix}1 & 0\\
k\epsilon & 1\end{smallmatrix}\right)$, then for any $\left(\begin{smallmatrix}t_{0} & 0\\
0 & t_{0}^{-1}\end{smallmatrix}\right)\in T$ we have\begin{eqnarray*}
L^{-1}(\left(\begin{smallmatrix}1 & 0\\
k\epsilon & 1\end{smallmatrix}\right)U_{2}) & \simeq & L^{-1}(F(\kappa)U_{2}F(\kappa)^{-1})=L^{-1}\left(F\big[\kappa\big(\begin{smallmatrix}t_{0} & 0\\
0 & t_{0}^{-1}\end{smallmatrix}\big)\big]U_{2}F\big[\kappa\big(\begin{smallmatrix}t_{0} & 0\\
0 & t_{0}^{-1}\end{smallmatrix}\big)\big]^{-1}\right)\\
 & \simeq & L^{-1}\left(L(\kappa\big(\begin{smallmatrix}t_{0} & 0\\
0 & t_{0}^{-1}\end{smallmatrix}\big))U_{2}\right)=L^{-1}\left(\big(\begin{smallmatrix}t_{0} & 0\\
0 & t_{0}^{-1}\end{smallmatrix}\big)\left(\begin{smallmatrix}1 & 0\\
k\epsilon & 1\end{smallmatrix}\right)F\big(\begin{smallmatrix}t_{0} & 0\\
0 & t_{0}^{-1}\end{smallmatrix}\big)U_{2}\right)\\
 & = & L^{-1}\left(\big(\begin{smallmatrix}t_{0} & 0\\
0 & t_{0}^{-1}\end{smallmatrix}\big)\left(\begin{smallmatrix}1 & 0\\
k\epsilon & 1\end{smallmatrix}\right)\big(\begin{smallmatrix}t_{0}^{q} & 0\\
0 & t_{0}^{-q}\end{smallmatrix}\big)U_{2}\right)=L^{-1}\left(\left(\begin{smallmatrix}t_{0}^{q-1} & 0\\
t_{0}^{q+1}k\epsilon & t_{0}^{1-q}\end{smallmatrix}\right)U_{2}\right).\end{eqnarray*}
Thus $L^{-1}(\left(\begin{smallmatrix}1 & 0\\
k\epsilon & 1\end{smallmatrix}\right)U_{2})\simeq L^{-1}(\big(\begin{smallmatrix}1 & 0\\
t_{0}^{2}k\epsilon & 1\end{smallmatrix}\big)U_{2})$, if $t_{0}^{q}=t_{0}$. We see from this that the equivalence class
of $L^{-1}(\left(\begin{smallmatrix}1 & 0\\
k\epsilon & 1\end{smallmatrix}\right)U_{2})$ does only depend on the coset of $k$ in $\mathbb{F}^{\times}/\mathbb{F}_{q}^{\times2}$.
In the following we will show the stronger result that all varieties
$X_{x}$ for $x=\left(\begin{smallmatrix}1 & 0\\
k\epsilon & 1\end{smallmatrix}\right)$, $k\in\mathbb{F}^{\times}$ afford the same irreducible representations.

In \cite{lusztigpreprint}, 3.3 it is claimed that the variety $X_{e}$
realises $2q$ irreducible representations of dimension $(q^{2}-1)/2$.
The calculations of Lusztig already show that this cannot be the case,
since the variety affords a permutation representation of dimension
$q(q^{2}-1)$ and any permutation representation contains at least
one copy of the trivial representation. We will now give the correct
decomposition of this permutation representation. Hence we will see
that most representations of dimension $(q^{2}-1)/2$ are not realised
in varieties of the form $X_{x}$, $x\in G_{2}$.

For any $k\in\mathbb{F}^{\times}$ we write $X_{k}=L^{-1}(\left(\begin{smallmatrix}1 & 0\\
k\epsilon & 1\end{smallmatrix}\right)U_{2})$, by abuse of notation. The variety $X_{k}$ is endowed with an action
of the group $U_{2}\cap\left(\begin{smallmatrix}1 & 0\\
k\epsilon & 1\end{smallmatrix}\right)U_{2}\left(\begin{smallmatrix}1 & 0\\
k\epsilon & 1\end{smallmatrix}\right)^{-1}=U_{2}^{1}$, acting by right multiplication. By Lemma \ref{lem:X/U} we have
$X_{k}\sim X_{k}/U_{2}^{1}$. \renewcommand{\theenumi}{\emph{\roman{enumi}}}

\begin{thm}
Let $k\in\mathbb{F}^{\times}$. The virtual $G_{2}^{F}$-representation
$\sum_{i\geq0}(-1)^{i}H_{c}^{i}(X_{k})$ decomposes into the direct
sum of the following representations
\begin{enumerate}
\item $4$ distinct irreducible representations of dimension $(q^{2}-1)/2$,
each one with multiplicity $(q-1)/2$, 
\item the irreducible representations of dimension $1,\  q,\ (q+1)/2$,
each with multiplicity one, 
\item the irreducible representations of dimension $q+1$, each with multiplicity
$2$.
\end{enumerate}
Moreover, for all $k,k'\in\mathbb{F}^{\times}$ we have $X_{k}\sim X_{k'}$.
\end{thm}
\begin{proof}
The proof goes as follows. First we calculate the variety $X_{k}$
explicitly following the calculations of Lusztig for the case $k=1$
(cf. \cite{lusztigpreprint}, 3.3), and show that the representations
afforded by $X_{k}$ can be realised by a finite ($0$-dimensional)
variety $\overline{X}_{k}$. Next we show that there exists a partition
of $\overline{X}_{k}$ into closed subset such that the action of
$G_{2}^{F}$ on the parts is transitive. We identify a part $\overline{X}_{k}^{(\pm1,0)}$
with stabiliser $S^{F}$, where $S^{F}$ is identical to the group
$\Stab_{G^{F}}(\chi_{\psi})\cdot N$ in Section \ref{sec:The-Reps-of}.
Now Lemma \ref{lem:partition} tells us that the $G_{2}^{F}$-representation
afforded by $\overline{X}_{k}$ is isomorphic to the representation
induced from the $S^{F}$-representation afforded by the part $\overline{X}_{k}^{(\pm1,0)}$.
Finally, we show how to decompose the latter representation with respect
to the action of a certain abelian group, and thanks to the results
of Section \ref{sec:The-Reps-of} we can identify exactly which representations
of $G_{2}^{F}$ occur in the cohomology of $X_{k}$.

Let $g=\begin{pmatrix}a & b\\
c & d\end{pmatrix}\in G_{2}$. The condition that $g\in X_{k}$ is that\[
F\begin{pmatrix}a & b\\
c & d\end{pmatrix}\in\begin{pmatrix}a & b\\
c & d\end{pmatrix}\begin{pmatrix}1 & 0\\
k\epsilon & 1\end{pmatrix}\begin{pmatrix}1 & x\\
0 & 1\end{pmatrix}=\begin{pmatrix}a+bk\epsilon & b+(a+bk\epsilon)x\\
c+dk\epsilon & d+(c+dk\epsilon)x\end{pmatrix},\]
for some $x\in\mathbb{F}[[\epsilon]]/(\epsilon^{2})$. This condition
is equivalent to the system of equations\[
\begin{cases}
F(a)=a+bk\epsilon,\  F(b)=b+(a+bk\epsilon)x,\\
F(c)=c+dk\epsilon,\  F(d)=d+(c+dk\epsilon)x.\end{cases}\]
In order to eliminate $x$ from the equations, we note that since
$g\in G_{2}$, the above system is equivalent to \[
\begin{cases}
F(a)=a+bk\epsilon,\  F(c)=c+dk\epsilon,\\
(F(b)-b)(c+dk\epsilon)=(F(d)-d)(a+bk\epsilon).\end{cases}\]
Setting $a=a_{0}+a_{1}\epsilon,b=b_{0}+b_{1}\epsilon,c=c_{0}+c_{1}\epsilon,d=d_{0}+d_{1}\epsilon$,
we obtain \[
\begin{cases}
a_{0}^{q}=a_{0},\  c_{0}^{q}=c_{0},\  a_{1}^{q}=a_{1}+b_{0}k,\  c_{1}^{q}=c_{1}+d_{0}k,\\
(b_{0}^{q}-b_{0})c_{0}=(d_{0}^{q}-d_{0})a_{0},\ (b_{1}^{q}-b_{1})(c_{1}+d_{0}k)=(d_{1}^{q}-d_{1})(a_{1}+b_{0}k).\end{cases}\]
Thus, we may identify $X_{k}$ with the set of all $(a_{0},b_{0},c_{0},d_{0},a_{1},b_{1},c_{1},d_{1})\in\mathbb{F}^{8}$
such that\[
\begin{array}{l}
(a)\  a_{0}^{q}=a_{0},\  c_{0}^{q}=c_{0},\  a_{1}^{q}=a_{1}+b_{0}k,\  c_{1}^{q}=c_{1}+d_{0}k,\\
(b)\  a_{0}d_{0}-b_{0}c_{0}=1,\  a_{0}d_{1}+a_{1}d_{0}-b_{0}c_{1}-b_{1}c_{0}=0\ (\Leftrightarrow\det g=1),\\
(c)\ (b_{0}^{q}-b_{0})c_{0}=(d_{0}^{q}-d_{0})a_{0},\\
\ \ \ \ \  b_{1}^{q}c_{0}-b_{1}c_{0}+b_{0}^{q}c_{1}-b_{0}c_{1}+b_{0}^{q}d_{0}k=d_{1}^{q}a_{0}+d_{1}^{q}a_{0}-a_{0}d_{1}+d_{0}^{q}a_{1}-d_{0}a_{1}+d_{0}^{q}b_{0}k.\end{array}\]
Now the first equation $(c)$ follows by raising the first equation
$(b)$ to the $q$th power, and using the two first equations $(a)$.
The second equation $(c)$ follows by raising the second equation
$(b)$ to the $q$th power, using the two last equations $(a)$, and
adding again the second equation $(b)$. Thus the equations $(c)$
can be omitted.

The first equation $(b)$ can be written (using $(a)$):\[
a_{0}(c_{1}^{q}-c_{1})k^{-1}-c_{0}(a_{1}^{q}-a_{1})k^{-1}=1,\textrm{ that is }(a_{0}c_{1}-c_{0}a_{1})^{q}-(a_{0}c_{1}-c_{0}a_{1})=k.\]
Setting $f=a_{0}c_{1}-c_{0}a_{1}$, we see that $X_{k}$ can be identified
with the set of all $(a_{0},b_{0},c_{0},d_{0},a_{1},b_{1},c_{1},d_{1},f)\in\mathbb{F}^{9}$
such that\[
\begin{cases}
a_{0}^{q}=a_{0},\  c_{0}^{q}=c_{0},\  a_{1}^{q}=a_{1}+b_{0}k,\  c_{1}^{q}=c_{1}+d_{0}k,\\
f^{q}-f=k,\  f=a_{0}c_{1}-c_{0}a_{1},\  a_{0}d_{1}+a_{1}d_{0}-b_{0}c_{1}-b_{1}c_{0}=0.\end{cases}\]
We now factor out by the action of $U_{2}^{1}$. If $u=\left(\begin{smallmatrix}1 & x\epsilon\\
0 & 1\end{smallmatrix}\right)\in U_{2}^{1}$, then the action $g\mapsto gu$ on $X_{k}$ is given in terms of
coordinates by\[
(a_{0},b_{0},c_{0},d_{0},a_{1},b_{1},c_{1},d_{1},f)\longmapsto(a_{0},b_{0},c_{0},d_{0},a_{1},b_{1}+a_{0}x,c_{1},d_{1}+c_{0}x,f).\]
Suppose that $(a_{0},b_{0},c_{0},d_{0},a_{1},b_{1},c_{1},d_{1},f)$
and $(a_{0},b_{0},c_{0},d_{0},a_{1},b_{1}',c_{1},d_{1}',f)$ are two
points on $X_{k}$. Then \[
\begin{array}{l}
a_{0}d_{1}+a_{1}d_{0}-b_{0}c_{1}-b_{1}c_{0}=0,\textrm{ and }a_{0}d_{1}'+a_{1}d_{0}-b_{0}c_{1}-b_{1}'c_{0}=0.\end{array}\]
These equations imply that $a_{0}(d_{1}-d_{1}')=c_{0}(b_{1}-b_{1}')$,
which is equivalent to $b_{1}'=b_{1}+a_{0}x$, $d_{1}'=d_{1}'+c_{0}x$,
for some $x\in\mathbb{F}$. Thus the quotient variety $X_{k}/U_{2}^{1}$
may be identified with the set of points $(a_{0},b_{0},c_{0},d_{0},a_{1},c_{1},f)\in\mathbb{F}^{7}$
such that \[
a_{0}^{q}=a_{0},\  c_{0}^{q}=c_{0},\  a_{1}^{q}=a_{1}+b_{0}k,\  c_{1}^{q}=c_{1}+d_{0}k,\  f^{q}-f=k,\  f=a_{0}c_{1}-c_{0}a_{1}.\]
This in turn is naturally isomorphic to\[
\{(a_{0},c_{0},a_{1},c_{1},f)\in\mathbb{F}^{5}\mid a_{0}^{q}=a_{0},\  c_{0}^{q}=c_{0},\  f^{q}-f=k,\  f=a_{0}c_{1}-c_{0}a_{1}\}.\]
We consider the obvious projection $\alpha:(a_{0},c_{0},a_{1},c_{1},f)\mapsto(a_{0},c_{0},f)$
of this set, to the finite set \[
\overline{X}_{k}=\{(a_{0},c_{0},f)\in\mathbb{F}^{3}\mid a_{0}^{q}=a_{0},\  c_{0}^{q}=c_{0},\  f^{q}-f=k,\ (a_{0},c_{0})\neq(0,0)\}.\]
We remark that if $f^{q}-f=k$ for some $f\in\mathbb{F}$, then $(f+f_{0})^{q}-(f+f_{0})=k$
for any $f_{0}\in\mathbb{F}_{q}$. Hence, for any $k\in\mathbb{F}^{\times}$,
the equation $f^{q}-f=k$ has $q$ solutions, and if $f$ is such
a solution, then every solution is of the form $f+f_{0}$, for some
$f_{0}\in\mathbb{F}_{q}$.

Define an action of $G_{2}^{F}$ on $\overline{X}_{k}$ as the unique
action such that $\alpha(gx)=g\alpha(x)$, for $x\in X_{k}$, $g\in G_{2}^{F}$.
The fibre of $\alpha$ at $(a_{0},c_{0},f)\in\overline{X}_{k}$ is
the affine line $\{(a_{1},c_{1})\in\mathbb{F}^{2}\mid a_{0}c_{1}-c_{0}a_{1}=f\}$.
Thus by Lemma \ref{lem:X/U}, $X_{k}\sim X_{k}/U_{2}^{1}\sim\overline{X}_{k}$.
Since $\overline{X}_{k}$ is a $0$-dimensional variety, it follows
from Lemma \ref{lem:fixed-point} that $\sum_{i\geq0}(-1)^{i}H_{c}^{i}(\overline{X}_{k})=H_{c}^{0}(\overline{X}_{k})\simeq\overline{\mathbb{Q}}_{l}[\overline{X}_{k}]$,
which is a permutation representation of dimension $|\overline{X}_{k}|=q(q^{2}-1)$.

We now turn to the problem of decomposing the representation $\overline{\mathbb{Q}}_{l}[\overline{X}_{k}]$
into irreducibles. Consider the group\[
A=\Bigl{\{}\begin{pmatrix}\pm1 & 0\\
x\epsilon & \pm1\end{pmatrix}\mid x\in\mathbb{F}_{q}\Bigr{\}}.\]
This group acts on $X_{k}$ by right multiplication. There is a unique
action of $A$ on $\overline{X}_{k}$ satisfying $\alpha(xa)=\alpha(x)a$,
for $x\in X_{k}$, $a\in A$. The action on $\overline{X}_{k}$ is
given in terms of coordinates by\[
(a_{0},c_{0},f)\longmapsto(\pm a_{0},\pm b_{0},f+x).\]
The set of orbits $\overline{X}_{k}/A$ defines a partition of $\overline{X}_{k}$
into closed subsets $\overline{X}_{k}^{(a_{0},c_{0})}$, indexed by
pairs $(a_{0},c_{0})\in\mathbb{F}_{q}^{2}/\{\pm1\}$, $(a_{0},c_{0})\neq(0,0)$.
Hence, each orbit contains $2q$ elements.

Now consider the action of the group $G_{2}^{F}$ on $\overline{X}_{k}$.
For $g=\left(\begin{smallmatrix}x_{0}+x_{1}\epsilon & y_{0}+y_{1}\epsilon\\
z_{0}+z_{1}\epsilon & w_{0}+w_{1}\epsilon\end{smallmatrix}\right)\in G_{2}^{F}$, the action is given in terms of coordinates by\[
\begin{array}{l}
(a_{0},c_{0},f)\longmapsto(x_{0}a_{0}+y_{0}c_{0},z_{0}a_{0}+w_{0}c_{0},\\
f+a_{0}^{2}(x_{0}z_{1}-z_{0}x_{1})+a_{0}c_{0}(x_{0}w_{1}+y_{0}z_{1}-z_{0}y_{1}-w_{0}x_{1})+c_{0}^{2}(y_{0}w_{1}-w_{0}y_{1})).\end{array}\]
Thus $g\overline{X}_{k}^{(a_{0},c_{0})}=\overline{X}_{k}^{(x_{0}a_{0}+y_{0}c_{0},z_{0}a_{0}+w_{0}c_{0})}$,
and $G_{2}^{F}$ acts transitively on the set of orbits $\overline{X}_{k}^{(a_{0},c_{0})}$.
The stabiliser of the orbit $\overline{X}_{k}^{(\pm1,0)}$ is given
by\[
S^{F}:=\Stab_{G_{2}^{F}}(\overline{X}_{k}^{(\pm1,0)})=\Bigl{\{}\begin{pmatrix}\pm1+x_{1}\epsilon & y_{0}+y_{1}\epsilon\\
z_{1}\epsilon & \pm1+w_{1}\epsilon\end{pmatrix}\in G_{2}^{F}\Bigr{\}}.\]
It follows from Lemma \ref{lem:partition} that the $G_{2}^{F}$-representation
$\sum_{i\geq0}(-1)^{i}H_{c}^{i}(X_{k})\simeq H_{c}^{0}(\overline{X}_{k})$
is induced by the $S^{F}$-representation $\sum_{i\geq0}(-1)^{i}H_{c}^{i}(\overline{X}_{k}^{(\pm1,0)})=H_{c}^{0}(\overline{X}_{k}^{(\pm1,0)})$.
We will determine the latter by using some character theory.

Let $s$ denote an element in $S^{F}$, and let $\chi:A\rightarrow\overline{\mathbb{Q}}_{l}^{\times}$
be a character. Since induction preserves direct sums, we have \[
H_{c}^{0}(\overline{X}_{k})=\bigoplus_{\chi}\Ind_{S^{F}}^{G_{2}^{F}}H_{c}^{0}(\overline{X}_{k}^{(\pm1,0)})_{\chi},\]
and by Lemma \ref{lem:char-formula} and Lemma \ref{lem:fixed-point}
we have \[
\Tr(s,H_{c}^{0}(\overline{X}_{k}^{(\pm1,0)})_{\chi})=\frac{1}{|A|}\sum_{a\in A}\chi(a^{-1})|(\overline{X}_{k}^{(\pm1,0)})^{(s,a)}|.\]
Thus we have to determine the fixed points of $\overline{X}_{k}^{(\pm1,0)}$
under the action of $(s,a)\in S^{F}\times A$. Every $s\in S^{F}$
can be decomposed uniquely as $s=ua'$, where $u\in\{\left(\begin{smallmatrix}1+x\epsilon & y_{0}+y_{1}\epsilon\\
0 & 1-x\epsilon\end{smallmatrix}\right)\in G_{2}^{F}\}$, and $a'\in A$. The set $\overline{X}_{k}^{(\pm1,0)}$ is fixed
under conjugation by elements in $A$, and the element $u$ leaves
it fixed. Thus\[
|(\overline{X}_{k}^{(\pm1,0)})^{(s,a)}|=\begin{cases}
2q & \textrm{if }a'=a^{-1},\\
0 & \textrm{otherwise},\end{cases}\]
and so\[
\Tr(s,H_{c}^{0}(\overline{X}_{k}^{(\pm1,0)})_{\chi})=\frac{1}{2q}\chi(a')2q=\chi(a').\]
We see that the character of the representation $H_{c}^{0}(\overline{X}_{k}^{(\pm1,0)})$
is the direct sum of all characters of $S^{F}$ given by\[
\begin{pmatrix}\pm1+x_{1}\epsilon & y_{0}+y_{1}\epsilon\\
z_{1}\epsilon & \pm1+w_{1}\epsilon\end{pmatrix}\longmapsto\chi^{\pm}(\pm1)\cdot\psi(z_{1}),\]
where $\chi^{\pm}\in\Hom(\{\pm1\},\overline{\mathbb{Q}}_{l}^{\times}),\ \psi\in\Hom(\mathbb{F}_{q}^{+},\overline{\mathbb{Q}}_{l}^{\times}).$
Comparing this with the description of the irreducible representations
of $G_{2}^{F}$ of dimension $(q^{2}-1)/2$ given in Section \ref{sec:The-Reps-of},
we see that there are four characters $\chi$ for which the induced
characters of $G_{2}^{F}$ are distinct and irreducible. All other
characters $\chi$ such that $\chi^{2}\neq1$ give rise to representations
of $G_{2}^{F}$ isomorphic to one of these. Next consider $\chi=1$.
Inducing the trivial character of $S^{F}$ to $G_{2}^{F}$ gives the
character of the permutation representation $\overline{\mathbb{Q}}_{l}[G_{2}^{F}/S^{F}]\simeq\overline{\mathbb{Q}}_{l}[G^{F}/\{\left(\begin{smallmatrix}\pm1 & y_{0}\\
0 & \pm1\end{smallmatrix}\right)\}]$, and from the classical finite field case we know that this representation
consists of the representations of dimension $1$ and $q$, and $(q-3)/2$
irreducible representations of dimension $q+1$. Analogously, inducing
the character $\chi$ for which $\chi^{2}=1,\ \chi\neq1$, we see
that the resulting representation consists of the two irreducible
representations of dimension $(q+1)/2$, and $(q-3)/2$ irreducible
representations of dimension $q+1$ isomorphic to the ones corresponding
to $\chi=1$.

This discussion shows that the above results are independent of the
choice of $k\in\mathbb{F}^{\times}$. Hence, the theorem is proved.
\end{proof}

\section{Realising The Missing Representations}

Using the observations of the previous sections, we show how all irreducible
representations of dimension $(q^{2}-1)/2$ of $G_{2}^{F}$ can be
realised in the cohomology of a certain variety.

Consider the subgroup $S=ZG_{2}^{1}U$ of $G_{2}$, where $Z$ is
the centre of $G_{2}$. Then it is clearly an $F$-stable subgroup,
and $S^{F}$ is compatible with the notation used in the previous
section. As we have seen, all irreducible representations of dimension
$(q^{2}-1)/2$ of $G_{2}^{F}$ can be obtained by inducing certain
$1$-dimensional representations of the subgroup $S^{F}$. Some elementary
calculations show that the commutator subgroup of $S$ is given by\[
(S,S)=\Bigl\{\begin{pmatrix}1+x\epsilon & y\epsilon\\
0 & 1-x\epsilon\end{pmatrix}\Bigr\},\]
and this is again an $F$-stable subgroup. Consider the variety\[
Y=\{ g\in G_{2}\mid g^{-1}F(g)\in(S,S)\}/(S,S),\]
where $(S,S)$ acts by right multiplication. The variety $Y$ is acted
on by $G_{2}^{F}$ and $S^{F}$ by left and right multiplication,
respectively.

Since $(S,S)$ is $F$-stable, there is an obvious inclusion map $Y\rightarrow G_{2}/(S,S)$,
whose image is $(G_{2}/(S,S))^{F}\simeq G_{2}^{F}/(S,S)^{F}$. We
thus have an isomorphism $Y\simeq G_{2}^{F}/(S,S)^{F}$. 

Let $\chi$ be a character of $S^{F}$ such that $\chi$ is an extension
of a character $\left(\begin{smallmatrix}1 & 0\\
x\epsilon & 1\end{smallmatrix}\right)\mapsto\psi(x)$, for some $\psi\in\Hom(\mathbb{F}_{q}^{+},\overline{\mathbb{Q}}_{l}^{\times})$.
We have seen that there are two distinct orbits of such characters
$\chi$, with respect to the action of $G^{F}$. Let $\theta$ be
an arbitrary character of the subgroup $\{\left(\begin{smallmatrix}\pm1 & y_{0}\\
0 & \pm1\end{smallmatrix}\right)\in G^{F}\}$. 

\begin{prop}
Every irreducible representation of $G_{2}^{F}$ of dimension $(q^{2}-1)/2$
is given by\[
\sum_{i\geq0}(-1)^{i}(H_{c}^{i}(Y)\otimes\chi\otimes\theta)=H_{c}^{0}(Y)\otimes\chi\otimes\theta,\]
where the characters $\chi$ and $\theta$ are identified with their
corresponding $1$-dimensional representations of $S^{F}$. Moreover,
if $\chi$ runs through a set of representatives of orbits under $G^{F}$,
then each irreducible representation of dimension $(q^{2}-1)/2$ appears
exactly once.
\end{prop}
\begin{proof}
We have seen above that $Y\simeq G_{2}^{F}/(S,S)^{F}$, so by Lemma
\ref{lem:fixed-point} the cohomology of $Y$ is concentrated in degree
$0$, and $H_{c}^{0}(Y)\simeq\overline{\mathbb{Q}}_{l}[G_{2}^{F}/(S,S)^{F}]$.
Now let $V$ be a representation of $S^{F}$ that factors through
$(S,S)^{F}$, i.e., a $1$-dimensional representation. Then there
is a representation of $G_{2}^{F}$ given by $H_{c}^{0}(Y)\otimes_{\overline{\mathbb{Q}}_{l}[S^{F}]}V$,
and this representation is isomorphic to the one given by inducing
$V$ to $G_{2}^{F}$ (cf. \cite{dignemichel}, ch.~4). By the description
of the irreducible representations of dimension $(q^{2}-1)/2$ as
induced by $1$-dimensional representations of $S^{F}$ given in Section
\ref{sec:The-Reps-of}, it follows that these representations are
given by $H_{c}^{0}(Y)\otimes\chi\otimes\theta$, as asserted.
\end{proof}
\bibliographystyle{alex}
\bibliography{alex}

\end{document}